\font\tenmsb=msbm10
   \font\sevenmsb=msbm7
   \font\fivemsb=msbm5

\newfam\msbfam
      \textfont\msbfam=\tenmsb
      \scriptfont\msbfam=\sevenmsb
      \scriptscriptfont\msbfam=\fivemsb
\def\Bbb#1{{\fam\msbfam #1}}

\font\teneufm=eufm10
\font\seveneufm=eufm7
\font\fiveeufm=eufm5
\newfam\eufmfam
\textfont\eufmfam=\teneufm
\scriptfont\eufmfam=\seveneufm
\scriptscriptfont\eufmfam=\fiveeufm
\def\frak#1{{\fam\eufmfam\relax#1}}

\newcommand\qed{{\hspace*{\fill}Q.E.D.\vskip12pt plus 1pt}}

\newcommand\sB{{\cal B}}
\newcommand\sD{{\cal D}}
\newcommand\sE{{\cal E}}
\newcommand\sF{{\cal F}}
\newcommand\sJ{{\cal J}}
\newcommand\sI{{\cal I}}
\newcommand\sK{{\cal K}}
\newcommand\sL{{\cal L}}
\newcommand\sO{{\cal O}}
\newcommand\sT{{\cal T}}
\newcommand\sX{{\cal X}}
\newcommand\sZ{{\cal Z}}

\newcommand\grk{\kappa}
\newcommand\grs{\sigma}
\newcommand\vphi{\varphi}

\newcommand\comp{{\Bbb C}}

\newcommand\hatX{{\widehat X}}
\newcommand\hatL{{\widehat L}}
\newcommand\pn[1]{{\Bbb P}^{#1}}
\newcommand\pnsheaf[2]{\sO_{{\Bbb P}^{#1}}\left({#2}\right)}
\newcommand\pnpair[2]{\left(\pn {#1},\pnsheaf {#1}{#2}\right)}
\newcommand\grass[2]{{\mbox{\rm Gr}\left(#1,#2\right)}}

\newcommand\proj[1]{{\Bbb P}({#1})}
\newcommand\proof{\noindent{\em Proof.}\ \ }
\newcommand\bx{{\:\:\!\!\sqsubset\!\;\!\!\;\!\!\!\!\!\sqsupset\!\!\!\:\:\!\!
\!\!\!\!\times}\, }

\def\red{{\mathop{\rm red}\nolimits}}
\def\rank{{\mathop{\rm rank}\nolimits}}

\documentstyle[12pt,twoside]{article}

\newtheorem{theorem}{{\rm T\sc heorem}}[section]
\newtheorem{lemma}[theorem]{{\rm L\sc emma}}
\newtheorem{corollary}[theorem]{{\rm C\sc orollary}}
\newtheorem{proposition}[theorem]{{\rm P\sc roposition}}

\newtheorem{deff}[theorem]{{\rm D\sc efinition}}
\newtheorem{re}[theorem]{{\rm R\sc emark}}
\newtheorem{pargrph}[theorem]{}

\newenvironment{definition}{\begin{deff}\em}{\end{deff}}
\newenvironment{rem*}{\begin{re}\em}{\end{re}}
\newenvironment{prgrph*}[1]{\indent\begin{pargrph}{\bf #1.}\em\
}{\end{pargrph}}

\begin{document}

\title{On higher order embeddings of Fano threefolds by the anticanonical
linear system}
\author{Mauro C. Beltrametti, Sandra Di Rocco, and Andrew
J. Sommese}

\maketitle

\begin{abstract}
In this article the map given by the anticanonical bundle of a
Fano manifold is investigated with respect to a number of natural notions 
of higher order embeddings of projective manifolds. This is of
importance in the understanding of higher order embeddings of the special
varieties of adjunction theory, which are usually fibered by special Fano
manifolds. An analysis is carried out of the higher order embeddings of the
special varieties of adjunction theory that arise in the study of the first
and second reductions. Special attention is given to determining the order
of the anticanonical embeddings of the   three dimensional Fano manifolds
which have been classified by Iskovskih, Mori, and Mukai and also of the
Fano complete intersections in $\pn N$.\\
\  \\
\noindent{1991 {\em Mathematics Subject Classification}. Primary 14J45, 14J40;
Secondary 14M10, 14N99.\newline
}
\end{abstract}

\section*{Introduction}
Let $X$ be an $n$-dimensional connected complex projective manifold.
 There are  three natural notions (see (\ref{kthemb})) of the
``order'' of an embedding given by a line bundle $L$.  From strongest to
 weakest they are
$k$-jet ampleness \cite{Plenum}, $k$-very ampleness \cite{BSCo}, and
$k$-spannedness \cite{Duke}.
For $k=1$ (respectively $k=0$) all three notions are equivalent to
very ampleness (respectively spannedness at all points of $X$ by global
sections).
  We consider the most natural notion to be $k$-very ample,
which by definition means that given any $0$-dimensional subscheme $\sZ$ of $X$
of  length $k+1$, the    map $ H^0(X,L)\to H^0(\sZ,L_{\sZ})$, is onto.
There has been considerable work on deciding the order of the embeddings
relative to these notions for the line bundles that come up on the standard
classes of varieties.

In this article we investigate what the ``order'' of the embedding
by $|-K_X|$ is.

In \S \ref{Background} we recall some known results that we will need in
the sequel.  We
also prove some foundational results on $k$-very ampleness that are not in the
literature, e.g., Lemma (\ref{Tensor}) and the useful lower bound for the
degree and the number of sections of $L$, Proposition (\ref{h0Est}).

In \S \ref{FCIsection} we  give
the $k$-jet ampleness, $k$-very ampleness,
and $k$-spannedness of $-K_X$ for Fano manifolds  $X\subset \pn N$
 which are   complete
intersections of hypersurfaces in $\pn N$.

In \S \ref{kAd}  we work out the order of the embeddings of the special
varieties
that come up in the study of the first and second reductions of
adjunction theory.  Fano manifolds of   special types come up naturally
as the fibers of degenerate adjunction morphisms.

In \S \ref{Mukai} we  continue the investigation of the Fano manifolds that came
up in \S \ref{kAd}, those with $-K_X\cong (n-2)L$ and $n\ge 4$.

In \S \ref{FanoS}, by using Fujita's classification of Del Pezzo threefolds,
 we completely settle the order of the embeddings of
$3$-dimensional Fano manifolds with very ample anticanonical bundle $-K_X$.

For further discussion and a guide to most of the published papers on $k$-very
ampleness, we refer to the book \cite{Book} by the first and
third author.  We also call attention to \cite{Sandra1} of the
second author,
which do a thorough investigation for surfaces of the questions analogous to
those we ask for
higher dimensions.

The third author would like to thank the University of Notre Dame and the
Alexander von Humboldt Stiftung for their support.

We would like to thank the referee for helpful comments.

\section{Background material}\label{Background}
\setcounter{theorem}{0}
Throughout this paper we deal with
complex projective manifolds $V$. We denote by $\sO_V$  the structure
sheaf of $V$ and by $K_V$ the canonical bundle. For any coherent sheaf
$\sF$ on $V$, $h^i(\sF)$ denotes the complex dimension of $H^i(V,\sF)$.

Let $L$ be a line bundle on $V$. $L$ is said to be {\em numerically
effective} ({\em nef}, for short) if $L\cdot C\geq 0$ for all effective
curves $C$ on $V$. $L$ is said to be {\em big} if $\grk(L)=\dim V$, where
$\grk(L)$ denotes the Kodaira dimension of $L$. If $L$ is nef then this
is equivalent to be $c_1(L)^n>0$, where $c_1(L)$ is the first Chern class
of $L$ and $n=\dim V$.
\begin{prgrph*}{Notation}\label{Notation}
In this paper, we use the standard
 notation  from algebraic geometry. Let us
only fix the following.
\begin{enumerate}
\item[] $\approx$ (respectively $\sim$) denotes  linear (respectively numerical)
equivalence of line bundles;
\item[] $|L|$, the complete linear system associated with a line bundle
$L$ on a variety $V$;
\item[] $\Gamma(L)=H^0(L)$ denotes the space of the global sections of $L$. We
say that $L$ is spanned if it is spanned at all points of $V$ by
$\Gamma(L)$;
\item[] $\grk(V):=\grk(K_V)$ denotes the Kodaira dimension, for $V$ smooth;
\item[] $b_2(V)=\sum_{p+q=2}h^{p,q}$ denotes the second Betti number of $V$, for
$V$ smooth, where $h^{p,q}:=h^q(\Omega_V^p)$ denotes the Hodge   $(p,q)$
number of $V$.
\end{enumerate}

Line bundles and divisors are used with little (or no) distinction. Hence
we   freely use the additive notation.
\end{prgrph*}
\begin{prgrph*}{Reductions}\label{Red} (See e.g., \cite[Chapters 7,
12]{Book}) Let $(\widehat{X},\widehat{L})$ be a smooth projective variety
of dimension $n\geq 2$ polarized with a very ample line bundle
$\widehat{L}$. A smooth polarized variety $(X,L)$ is called a  {\em
reduction} of
 $(\widehat{X},\widehat{L})$ if there is a
morphism $r:\widehat{X}\to X$ expressing $\widehat{X}$ as the
blowing up of $X$ at a finite set of points, $B$, such that
  $L:=(r_*\widehat{L})^{**}$ is an ample line bundle and $\widehat{L}\approx
r^*L-[r^{-1}(B)]$ (or, equivalently,
$K_{\widehat{X}}+(n-1)\widehat{L}\approx r^*(K_X+(n-1)L)$).

Note that there is a one to one correspondence between smooth divisors of
$|L|$ which contain the set $B$ and smooth divisors of $|\widehat{L}|$.

Except for an explicit list of well understood pairs $(\widehat{X},\widehat{L})$
(see in particular \cite[\S\S 7.2, 7.3]{Book}) we can
assume:
\begin{enumerate}
\item[{\rm a)}] $K_{\widehat{X}}+(n-1)\widehat{L}$ is spanned and big, and
$K_X+(n-1)L$ is very ample. Note that in this case this reduction, $(X,L)$,
is unique up to isomorphism. We will refer to it as {\em the first
reduction} of  $(\widehat{X},\widehat{L})$.
\item[{\rm b)}] $K_X+(n-2)L$ is nef and big, for $n\geq 3$.
\end{enumerate}

Then from the Kawamata-Shokurov basepoint free theorem (see e.g.,
\cite[(1.5.2)]{Book}) we know that $|m(K_X+(n-2)L|$, for $m\gg 0$, gives rise
to a morphism $\vphi:X\to Z$, with connected fibers and normal image. Thus
there is an ample line bundle $\sK$ on $Z$ such that $K_X+(n-2)L\approx
\vphi^*\sK$. Let $\sD:=(\vphi_*L)^{**}$. The pair $(Z,\sD)$, together with
the morphism $\vphi:X\to Z$ is called the {\em second reduction} of
$(\widehat{X},\widehat{L})$. The morphism $\vphi$ is very well behaved (see
e.g., \cite[\S\S 7.5, 7.6 and \S\S 12.1, 12.2]{Book}). In particular $Z$ has
terminal, $2$-factorial isolated singularities and $\sK\approx K_Z+(n-2)\sD$.
Moreover $\sD$ is a  $2$-Cartier divisor such that $2L\approx
\vphi^*(2\sD)-\Delta$, for some effective Cartier divisor $\Delta$ on $X$
which is $\vphi$-exceptional (see \cite[(7.5.6), (7.5.8)]{Book}).
\end{prgrph*}
\begin{prgrph*}{Nefvalue} (See e.g., \cite[\S 1.5]{Book}) Let $V$ be a smooth
projective variety and let $L$ be an ample line bundle on $V$. Assume that
$K_V$ is not nef. Then from the Kawamata rationality theorem (see e.g.,
\cite[(1.5.2)]{Book}) we know that there exists a rational number $\tau$ such
that $K_V+\tau L$ is nef and not ample. Such a number, $\tau$, is called the
{\em nefvalue} of $(V,L)$.

Since $K_V+\tau L$ is nef, it follows from the Kawamata-Shokurov base point
free theorem (see e.g., \cite[(1.5.1)]{Book}) that $|m(vK_V+uL)|$ is
basepoint free for all $m\gg 0$, where $\tau=u/v$. Therefore, for such $m$,
$|m(K_V+\tau L)|$ defines a morphism $f:V\to {\Bbb P}_{\comp}$. Let $f=s\circ
\Phi$ be the Remmert-Stein factorization of $f$, where $\Phi:V\to Y$ is a
morphism with connected fibers onto a normal projective variety $Y$ and
$s:Y\to  {\Bbb P}_{\comp}$ is a finite-to-one morphism. For $m$ large enough
such a morphism, $\Phi$, only depends on $(V,L)$ (see \cite[\S 1.5]{Book}).
We call  $\Phi:V\to Y$ the {\em nefvalue morphism} of $(V,L)$.
\end{prgrph*}
\begin{prgrph*}{Special varieties}\label{Special} (See e.g., \cite[\S
3.3]{Book}) Let $V$
be a smooth variety of dimension $n$ and let $L$ be an ample line bundle on $V$.

We say that $V$ is a {\em Fano manifold} if $-K_V$ is ample. We say that $V$
is a {\em Fano
manifold of index } $i$ if $i$ is the largest positive integer such that
$K_V\approx -iH$
for some ample line bundle $H$ on $V$. Note that $i\leq n+1$ (see e.g.,
\cite[(3.3.2)]{Book}) and $n-i+1$ is referred to as the {\em co-index} of $V$.

We say that a Fano manifold, $(V,L)$, is a {\em Del Pezzo manifold}
(respectively a {\em Mukai
manifold}) if
$K_V\approx -(n-1)L$ (respectively $K_V\approx -(n-2)L$).

We also say that $(V,L)$ is a {\em scroll} (respectively a {\em quadric
fibration};
respectively a {\em Del Pezzo fibration}; respectively a {\em Mukai
fibration}) over a
normal variety $Y$ of dimension $m$ if there exists a surjective morphism
with connected
fibers $p:V\to Y$ such that $K_V+(n-m+1)L\approx p^*\sL$ (respectively
$K_V+(n-m)L\approx p^*\sL$; respectively   $K_V+(n-m-1)L\approx p^*\sL$;
respectively   $K_V+(n-m-2)L\approx p^*\sL$) for some ample line bundle
$\sL$ on $Y$.
\end{prgrph*}
\begin{prgrph*}{$k$-th order embeddings}\label{kthemb} Let $V$ be a smooth
algebraic
variety. We denote the Hilbert scheme of $0$-dimensional subschemes 
$(\sZ,\sO_{\sZ})$ of
$V$ with ${\rm length}(\sO_\sZ)=r$ by $V^{[r]}$. Since we are working in
characteristic zero, we have
${\rm length}(\sO_\sZ)=h^0(\sO_\sZ)$.

We say that a line bundle $L$ on $V$ is $k$-{\em very ample} if the
restriction map
$\Gamma(L)\to\Gamma(\sO_\sZ(L))$ is onto for any $\sZ\in V^{[k+1]}$. Note
that $L$ is
$0$-very ample if and only if $L$ is spanned by global sections, and $L$ is
$1$-very ample
if and only if $L$ is very ample. Note also that for smooth surfaces
with $k\leq 2$, $L$ being $k$-very ample is equivalent to $L$ being
$k$-{\em spanned} in
the sense of \cite{Duke}, i.e., $\Gamma(L)$ surjects on $\Gamma
(\sO_\sZ(L))$ for any {\em
curvilinear} $0$-cycle $\sZ \in V^{[k+1]}$, i.e., any $0$-dimensional
subscheme,
$\sZ\subset V$, such that ${\rm length}(\sO_\sZ)=k+1$ and $\sZ\subset C$
for some smooth
curve $C$ on $V$ (\cite[(0.4), (3.1)]{Duke}).

Let $x_1,\ldots,x_r$ be $r$ distinct points on $V$.
Let ${\frak m}_i$ be the maximal ideal sheaves of the points $x_i\in V$,
$i=1,\ldots,r$. Note that the stalk of ${\frak m}_i$  at $x_i$ is
 nothing but the maximal ideal, ${\frak m}_i\sO_{V,x_i}$,
of the local ring $\sO_{V,x_i}$, $i=1,\ldots,r$. Consider the $0$-cycle
$\sZ= x_1+\cdots+x_r$. We say that $L$ is $k$-{\em jet ample at} $\sZ$
if, for every $r$-tuple $(k_1,\ldots,k_r)$ of positive integers such that
$\sum_{i=1}^rk_i=k+1$, the restriction map
$$\Gamma(L)\to\Gamma(L\otimes (\sO_V/\otimes_{i=1}^r{\frak m}_i^{k_i}))
\left(\cong \oplus_{i=1}^r\Gamma(L\otimes (\sO_V/{\frak m}_i^{k_i})\right)$$
is onto. Here  ${\frak m}_i^{k_i}$ denotes
the $k_i$-th tensor power of ${\frak m}_i$.

We say that  $L$ is $k$-{\em jet ample} if, for any
$r\geq 1$ and any $0$-cycle $\sZ= x_1+\cdots+x_r$, where $x_1,\ldots,x_r$
are $r$
distinct points on $V$, the line bundle $L$ is $k$-jet ample at $\sZ$.

Note that $L$ is $0$-jet ample if and only if $L$ is spanned by its global
sections and $L$ is $1$-jet ample if and only if $L$ is very ample.

Note also that if $L$ is $k$-jet ample, then $L$ is $k$-very ample (see
\cite[(2.2)]{Plenum} and compare also with (\ref{SpecialProp})).

We will use over and over through the paper the fact
\cite[(1.3)]{BSMZ}, that if $L$ is a $k$-very ample line bundle on  $V$, then
$L\cdot C\geq k$ for each irreducible curve $C$ on $V$.

We refer to \cite{BSAq}, \cite{BSCo} and \cite{BSMZ},  and  \cite{Plenum} for
more on $k$-spannedness, $k$-very ampleness and $k$-jet ampleness
respectively.
\end{prgrph*}

\begin{definition}\label{image}  Let  $p:X\to Y$   be a holomorphic map
between complex projective schemes. Let  $\sZ$   be a  $0$-dimensional
subscheme of
$X$  defined by the ideal sheaf  $\sJ_\sZ$.   Then the  {\em
image} $p(\sZ)$  of  $\sZ$   is the
$0$-dimensional subscheme of  $Y$  whose defining ideal is
$\sI=\{g\in\sO_Y\;|\;g\circ p\in \sJ_\sZ\}$.
\end{definition}

We need the following general fact.

\begin{lemma}\label{GenFa} Let $p:X\to Y$ be a morphism of quasiprojective
varieties $X$,
$Y$. Let $\sZ$ be a  $0$-dimensional subscheme of $X$ of length $k$. Then
$p(\sZ)$ has length $\leq k$.
\end{lemma}
\proof Let $x_1,\ldots,x_t$ be $t$ distinct points such that ${\rm
Supp}(\sZ)=\{x_1,\ldots,x_t\}$. Let ${\rm length}(\sO_\sZ)=k$ and ${\rm
length}(\sO_{\sZ,x_i})=k_i$, where $\sO_{\sZ,x_i}$ denotes the stalk of
$\sO_\sZ$ at
$x_i$, $i=1,\ldots,t$. Then $k=\sum_{i=1}^tk_i$. Set $\sZ':=p(\sZ)$ and
let $\sJ_\sZ$,
$\sJ_{\sZ'}$ be the ideal sheaves of $\sZ$, $\sZ'$ respectively.

Arguing by contradiction, assume that ${\rm length}(\sO_{\sZ'})> k$. Then
there are $k+1$ linearly independent functions, $g_0=1,g_1,\ldots,g_k$, in
$\sO_{\sZ'}$.  Let $y_1,\ldots,y_s$, $s\leq t$, be the images of the points
$x_i$, $i=1,\ldots,t$. For each $j=1,\ldots,s$, consider the vector subspace
of $\comp^{k+1}$ defined by $$V_j:=\{(\lambda_0,\ldots,\lambda_k)\in
\comp^{k+1}\;|\;\sum_{i=0}^k\lambda_i g_i\in {\frak m}_j\},$$  where ${\frak
m}_j$ denote the ideal sheaf of $y_j$. Since
$$\dim V_j\geq k+1-k_j,\;j=1,\ldots,s,\ \ {\rm and}\ \ \sum_{j=1}^sk_j\leq
k<k+1$$
we conclude that
$$\dim\cap_{j=1}^s V_j\geq \sum_{j=1}^s\dim V_j-(s-1)(k+1)\ge
k+1-\sum_{j=1}^sk_j\ge 1.$$
Thus, there exist $\lambda_0,\ldots,\lambda_k\in \comp$ such that
$\sum_{i=0}^k\lambda_i
g_i=0$ at $y_j$ for each $j=1,\ldots,s$. Thus
$$p^*\left(\sum_{i=0}^k\lambda_i g_i\right)=\sum_{i=0}^k\lambda_i p^*g_i=0$$
at $x_i$ for each $i=1,\ldots,t$. It follows that $\sum_{i=0}^k\lambda_i
p^*g_i\in\sJ_\sZ$
and hence $\sum_{i=0}^k\lambda_i g_i\in\sJ_{\sZ'}$, or
$\sum_{i=0}^k\lambda_i g_i=0$ in
$\sO_{\sZ'}$. This contradicts the assumption that $1,g_1,\ldots,g_k$ are
linearly independent.\qed

If $X_1$, $X_2$ are projective schemes and $\sF_1$, $\sF_2$ are sheaves
on $X_1$, $X_2$
respectively, we will denote
$$\sF_1\bx\sF_2:=p_1^*\sF_1\otimes p_2^*\sF_2,$$
where $p_1$, $p_2$ are the projections on the two factors.

The following is the $k$-very ample version of Lemma (3.2) of
\cite{Embedding}.

\begin{lemma}\label{Tensor} Let $X_1$, $X_2$ be complex projective schemes
and $L_1$,
$L_2$ line bundles on $X_1$, $X_2$ respectively. For $i=1,2$ assume that
$L_i$ is
$k_i$-very ample and let $k:=\min\{k_1,k_2\}$. Then $L_1\bx L_2$ is
$k$-very ample on
$X_1\bx X_2$. Furthermore if $L_1\bx L_2$ is $k'$-very ample on $X_1\bx X_2$
then $L_1$, $L_2$ are $k'$-very ample. \end{lemma}
\proof Let $(\sZ,\sO_\sZ)$ be a $0$-dimensional subscheme of length $k+1$ on
$X_1\times X_2$.  Let $p_i:X_1\times X_2\to X_i$, $i=1,2$, the projections
on the two factors. Let $\sZ_i:={p_i}(\sZ)$ be the $0$-dimensional
subschemes of $X_i$ obtained as   images of $\sZ$, as in
(\ref{image}), and let $\sI_{{\sZ}_i}$ be the ideal sheaves
defining $\sZ_i$, $i=1,2$. Let $J_i:=p_i^*{\sI_{{\sZ}_i}}$, $i=1,2$. There
exist
generating  function germs $f\in J_i$  of type $f=g\circ p_i$, $g\in
\sI_{{\sZ}_i}$, $i=1,2$. Therefore for each point $z\in \sZ_\red$,
$f(z)=g(p_i(z))=0$, $i=1,2$. This means that $z$ belongs to the subscheme
of $X_1\times X_2$ defined by the ideal sheaf of the image of
$J_1\bx J_2$ in $\sO_{X_1\times X_2}$.
Hence we have an inclusion of ideal sheaves $(J_1,J_2)\subset \sI_{\sZ}$.
But $(J_1,J_2)$ defines the subscheme $\sJ$ whose structural sheaf is
$$\sO_{\sJ}=p_1^*\sO_{\sZ_1}\otimes
p_2^*\sO_{\sZ_2}=\sO_{\sZ_1}\bx\sO_{\sZ_2}.$$
Thus we have a surjection
\begin{equation}\label{surj}
\sO_{\sZ_1}\bx\sO_{\sZ_2}\to \sO_\sZ\to 0.
\end{equation}
On the other hand, by the Kunneth formula, we have
\begin{equation}\label{K1}
H^0(L_1\bx L_2) = H^0(p_1^*L_1\otimes p_2^* L_2)\\
=H^0(L_1)\otimes H^0(L_2),
\end{equation}
as well as,
\begin{eqnarray*}\label{K2}
H^0(L_1\bx L_2\otimes\sO_{\sZ_1}\bx\sO_{\sZ_2})&=&
H^0(\sO_{\sZ_1}(L_1)\bx\sO_{\sZ_2}(L_2))\\
&=&H^0(\sO_{\sZ_1}(L_1))\otimes H^0(\sO_{\sZ_2}(L_2)).
\end{eqnarray*}
Therefore, from (\ref{surj}) and (\ref{K1}), we get a surjection
\begin{equation}\label{s1}
H^0(\sO_{\sZ_1}(L_1))\otimes H^0(\sO_{\sZ_2}(L_2))\to
H^0(\sO_\sZ(L_1\bx L_2)).
\end{equation}
By Lemma (\ref{GenFa}) we have that $\sZ_1$ is of length $\leq k+1\leq
k_1+1$. Therefore,
since $L_1$ is $k_1$-very ample, the restriction map $H^0(L_1)\to
H^0(\sO_{\sZ_1}(L_1))$
is onto. Similarly we have that $H^0(L_2)\to H^0(\sO_{\sZ_2}(L_2))$ is
onto. Thus by using
(\ref{K1}) and (\ref{s1}) we get a surjection
$$H^0(L_1\bx L_2)\to H^0(\sO_\sZ(L_1\bx L_2))\to 0.$$
This shows that $L_1\bx L_2$ is $k$-very ample.

To show the last part of the statement note that ${L_1\bx L_2}_{|X_1\times
x_2}\cong L_1$ for each $x_2\in X_2$. Then $L_1$ is $k'$-very ample if $L_1\bx
L_2$ is $k'$-very ample. Similarly for $L_2$.
 \qed

The following result is a useful partial generalization of \cite[Lemma 1.1]
{BSMZ} (compare also with \cite[(3.1)]{Plenum}).

\begin{lemma}\label{kBlowupLemma}Let $L$ denote a $k$-very ample line
bundle on an
$n$-dimensional projective manifold $X$. Assume that $k\ge 2$ and
let $\pi : \hatX\to X$
denote the blowup of $X$ at a finite set  $\{x_1,\ldots,x_{k-1}\}\subset X$
of distinct
points.  Let $E_i:= \pi^{-1}(x_i)$ for $1\le i\le k-1$.
 Then $\pi^*L-(E_1+\cdots+E_{k-1})$ is very ample.
\end{lemma}
\proof Assume that $h^0(L)=N+1$, and we use $|L|$ to embed
$X$ into $\pn N$. Consider a linear subspace $\pn t\subset\pn N$, where $t\leq
k-1$. Assume that $\pn t$ met $X$ in a positive dimensional set. Then we
can assume without loss of generality that it is a curve. If not we can choose
a $\pn {t-1}$ contained in the $\pn t$ which will meet $X$ in a set of dimension
at most one less. Clearly $t>1$ since otherwise $\pn t$ would be a line
contained in $X$, contradicting the $k$-very ampleness assumption  with $k\geq
2$. Choose a $\pn {t-1}$ in the $\pn t$ meeting $X$ in a finite set. This is
possible since $\pn t$ meets $X$ in some points and a curve, say $C$. By using
the very ampleness assumption we conclude that the $\pn {t-1}$ meets the curve
$C$ in at most $t$ points. Choose a hyperplane $H\subset\pn N$ meeting the
$\pn t$ in the $\pn {t-1}$. Then $H$ meets the curve $C$ in at most $t\leq
k-1$ points. Since $H$ restricts to $L$ on $X$, this contradicts the fact that
$L\cdot C\geq k$.
Therefore we conclude that
 any $\pn t\subset \pn N$
with $t\le k-1$ meets $X$ scheme theoretically in  a $0$-dimensional
subscheme, say $\sX$. Furthermore  ${\rm length}(\sO_{\sX})\leq t+1$. Indeed
otherwise $X\cap \pn t$ would contain a $0$-cycle of length $t+2\leq k+1$
which spans a $\pn {t+1}$ since $L$ is $k$-very ample.

 Thus by taking the projective space
$P:=\pn {k-2}$ generated by  $\{x_1,\ldots,x_{k-1}\}$, we see that this  $\pn
{k-2}$ meets $X$ in precisely $\{x_1,\ldots,x_{k-1}\}$.  Thus the blowup  $\grs
: Z\to \pn N$  of $\pn N$ at this  $P$ has $\hatX$ as the proper transform of
$X$ and the induced map from $\hatX$ to $X$ is $\pi$.  Since $\grs^*\sO_{\pn
N}(1)-\grs^{-1}(P)$ is spanned by global sections and
$\hatL:=\pi^*L-(E_1+\cdots+E_{k-1})$ is  the pullback to $\hatX$ of
$\grs^*\pnsheaf N 1-\grs^{-1}(P)$ we see that $\hatL$ is spanned by global
sections.  Thus  we conclude  that global sections of $\hatL$ give a map
$\phi :\hatX\to \pn{N-k+1}$.

Applying the fact that the $\pn{k-1}$ generated by the image in $\pn N$ of any
 $0$-dimensional subscheme $\sZ\subset X$ of length $k$ meets $X$ scheme
theoretically in $\sZ$,
it follows that  the rational map from $\hatX$ to $\pn {N-k+1}$ induced by the
projection of $\pn N$ from $P$    is one-to-one on $\hatX$.

To see that  $\hatL$ is very ample we consider the points
$x\in \hatX\setminus\cup_{i=1}^{k-1} E_i$ and the points $x\in
\cup_{i=1}^{k-1}E_i$ separately.

First assume that $x\in\hatX\setminus\cup_{i=1}^{k-1} E_i$.  Choose a tangent
vector  $\tau_x$ at $x$. Let  $u_1,\ldots,u_n$ be a choice of local
coordinates  defined in a neighborhood of $x$, all zero at $x$, and with
$\tau_x$ tangent
 to the $u_n$ axis.
 Let ${\frak a}_x$ denote the ideal sheaf which equals $\sO_X$ away from $x$,
and  at $x$ is defined by $(u_1,\ldots,u_{n-1},u_n^2)$.
 We can choose a zero dimensional subscheme $\sZ$ of $X$ of length $k+1$ that
is  defined by
  ${\frak a}_x\otimes {\frak m}_{x_1}\otimes\cdots\otimes {\frak m}_{x_{k-1}}$.
Since $H^0(L)\to H^0(L\otimes \sO_{\sZ})$ is onto by the definition of $k$-very
ampleness,
we conclude  that the map $\phi$ given by   global sections of $\hatL$
has nonzero  differential
evaluated on the tangent vector  $\tau_x$. It follows that the global sections
of $\hatL$
 embed away from $\cup_{i=1}^{k-1}E_i$.

To finish consider a point   $x\in \cup_{i=1}^{k-1}E_i$. Without loss of
generality we can  assume, by relabeling if necessary, that $x\in E_1$. We thus
have $x_1=\pi(x)$.
   Choose a tangent vector
$\tau_x$ at $x$. Note that since  the line bundle $\hatL$ is spanned and since
$\hatL_{E_1}\cong\pnsheaf {n-1}1$,   the restriction
$\phi_{E_1}$ is an embedding. Thus we can assume that $\tau_x$ is not tangent
to $E_1$. Let
 $\tau=d\pi_x(\tau_x)\in \sT_{X,x_1}$, where $d\pi_x:\sT_{\hatX,x}\to
\sT_{X,x_1}$ is the differential map and $\sT_{\hatX,x}$, $\sT_{X,x_1}$ are the
tangent bundles to $\hatX$, $X$ respectively.  Let
$u_1,\ldots,u_n$ be a choice of local coordinates  defined in a neighborhood of
$x_1$, all zero at $x_1$,   with $\tau$ tangent
 to the $u_n$ axis, and such that the proper transform of the
$u_n$ axis is tangent to $\tau_x$.
 Let ${\frak b}_x$ denote the ideal sheaf which equals $\sO_X$ away from $x$,
and  at $x$ is defined by $(u_1,\ldots,u_{n-1},u_n^3)$.
 We can choose a zero dimensional subscheme $\sZ$ of $X$ of length $k+1$ that
is
defined by ${\frak b}_x$ if $k=2$ and by ${\frak b}_x\otimes {\frak
m}_{x_2}\otimes\cdots\otimes {\frak m}_{x_{k-1}}$ if
$k\ge 3$. Since $H^0(L)\to H^0(L\otimes \sO_\sZ)$ is onto by the definition of
$k$-very  ampleness,
we conclude   that global sections of $\hatL$ give a map with nonzero
differential
evaluated on the tangent vector  $\tau_x$. It follows that the global
 sections of $\hatL$
 embed $\hatX$.
\qed
\begin{rem*}There is a nice interpretation of this result in terms of $X^{[k]}$,
 the Hilbert scheme of  $0$-dimensional subschemes of $X$ of
length $k$. There is a natural line bundle $\sL$ on
$X^{[k]}$ induced by $L$ on $X$.  The fact that $L$
is $k$-very ample is equivalent to
$\sL$ being very ample \cite{CG}.  Given a set $\{x_1,\ldots,x_{k-1}\}\subset X$
of $k-1$ distinct points of $X$,
the subscheme of all $0$-dimensional subschemes of $X$ of length $k$ containing
$\{x_1,\ldots,x_{k-1}\}$ is naturally isomorphic to $X$ blown up at
$\{x_1,\ldots,x_{k-1}\}\subset X$.  Under this identification
the very ample line bundle $\sL$  restricts to the line bundle
$\pi^*L-(E_1+\cdots+E_{k-1})$ of the above lemma.
\end{rem*}

We have the following estimate for the degree and the number of sections of
$L$.
\begin{proposition}\label{h0Est}Let $L$ denote a $k$-very ample line bundle
on an $n$-dimensional projective manifold $X$. Assume that $k\ge 2$. Then
$L^n\ge 2^n+k-2$ and $h^0(L)\ge 2n +k-1$.  Moreover  $h^0(L)= 2n
+k-1$ implies that $L^n = 2^n+k-2$.
\end{proposition}
\proof Let $\pi : \hatX\to X$ denote the blowup of $X$ at a finite set
  $\{x_1,\ldots,x_{k-1}\}\subset X$ of distinct
points.   Let $E_i:= \pi^{-1}(x_i)$ for $1\le i\le k-1$.
  By Lemma (\ref{kBlowupLemma}),   $\pi^*L-(E_1+\cdots+E_{k-1})$  is
very ample, and thus $\pi^*L-(E_1+\cdots+E_{k-1})-E_1$ is spanned.
This implies $(\pi^*L-(2E_1+E_2\cdots+E_{k-1}))^n\ge 0$ and
thus that $L^n \geq 2^n+k-2$ since $E_i^n=(-1)^n$, $i=1,\ldots,k-1$. If
$h^0(L)\le 2n +k-1$, then since $h^0(\hatL)=h^0(L)-k+1$, we have
$h^0(\hatL)\le 2n$.
The argument in \cite[Theorem (4.4.1), i)]{BSS2} shows that
$h^0(\hatL)= 2n$, so that $h^0(L)= 2n+k-1$, and ${\hatL}^n=2^n-1$.  Since
${\hatL}^n=L^n-k+1$, we are done.
\qed

\section{Fano complete intersections}\label{FCIsection}\setcounter{theorem}{0}
Our first result gives
the $k$-jet ampleness, $k$-very ampleness,
and $k$-spannedness of $-K_X$ for Fano  manifolds  $X\subset \pn N$
 which are   scheme theoretically complete  intersections of hypersurfaces in
$\pn N$. We say that a curve $\ell$ on $X$ is a line if $\sO_X(1)\cdot \ell=1$.

\begin{theorem}\label{CompIntFano}Let
$X$ be a positive dimensional
connected projective submanifold of $\pn N$, which is a
complete intersection of hypersurfaces of $\pn N$ of degree $d_i$, $i=1,\ldots,
r:=N-\dim X$.   If   the anticanonical bundle, $-K_X$ is ample and if $X$ is
not a degree $2$ curve, then $X$ contains a line. In particular
$-K_X$ is $(N+1-\sum_{i=1}^rd_i)$-jet ample, but not
$(N+2-\sum_{i=1}^rd_i)$-spanned.
\end{theorem}
\proof Since the curve case of this result is trivial, assume that $\dim
X\ge 2$.
Since $X$ is a complete intersection,  $K_X =\pnsheaf N{-N-1+\sum_{i=1}^rd_i}$,
 where $d_1,\ldots,d_r$ are the degrees of the hypersurfaces which intersect
 transversely in $X$.
Since $-K_X$ is very ample, we conclude that $\sum_{i=1}^rd_i\le N$, and
$-K_X$ is $(N+1-\sum_{i=1}^rd_i)$-jet ample (see
\cite[Corollary (2.1)]{Plenum}).  If we show that $X$ contains  a line $\ell$
it will follow that $-K_X\cdot \ell=N+1-\sum_{i=1}^rd_i$, and thus $-K_X$
is not
$(N+2-\sum_{i=1}^rd_i)$-spanned.  Thus we must only show that $X$ contains a
line.

Let $G$ denote the Grassmannian $\grass 2 {N+1}$ of $2$-dimensional
complex
 vector subspaces of $\comp^{N+1}$.  Let $\sF$ denote the tautological
rank $2$ quotient bundle of $G\times \comp^{N+1}$. Note that $G\times
\comp^{N+1}$ is naturally identified with $G\times H^0(\pnsheaf N 1)$.  Under
this identification
$\proj \sF\subset G\times \pn N$ is identified with the universal
family of linear
$\pn 1$'s contained in $\pn N$.
The image in $s'\in H^0(\sF)$ of a section $s$ of $\pnsheaf N 1$
vanishes at points of $G$
corresponding to lines in $s^{-1}(0)$.  Further a section $s$ of
$\pnsheaf N d$ maps
naturally to a section $s'$ of $\sF^{(d)}$, the $d$-th symmetric
tensor product $\sF$.
Here $s'$ vanishes at points of $G$ corresponding to the lines
contained in $s^{-1}(0)$.
Thus if $X$ is defined  by sections $s_1,\ldots,s_r$ of
$\pnsheaf N {d_1},\ldots,\pnsheaf N {d_r}$, the lines on $X$
correspond to the common
zeroes of the images, $s'_1,\ldots,s'_r$, in $\sF^{(d_1)},\ldots,\sF^{(d_r)}$.
Thus if we
show  that
$$c_{\rank\sF^{(d_1)}}(\sF^{(d_1)})\wedge\cdots\wedge
c_{\rank\sF^{(d_r)}}(\sF^{(d_r)})$$
is a nontrivial cohomology class, then it follows that
 $s'_1,\ldots,s'_r$ must have
common zeroes and $X$ must contain lines.

Note that $\rank\sF^{(d_i)}=d_i+1$. For odd $d_i$ we have
\begin{equation}\label{OddchernClasses}
c_{d_i+1}(\sF^{(d_i)})=
(d_i+1)^2c_2(\sF)\prod_{t=1}^{\frac{d_i-1}{2}}
\left(t(d_i-t)c_1^2(\sF)+(d_i-2t)^2c_2(\sF)\right)
\end{equation}
and for even $d_i$ we have
\begin{equation}\label{EvenchernClasses}
c_{d_i+1}(\sF^{(d_i)})=(d_i+1)^2c_2(\sF)
\frac{d_i}{2}c_1(\sF)\prod_{t=1}^{\frac{d_i}{2}-1}
\left(t(d_i-t)c_1^2(\sF)+(d_i-2t)^2c_2(\sF)\right).
\end{equation}
Since $\sF$ is spanned both $c_1(\sF)$ and $c_2(\sF)$ are semipositive classes
(see \cite[Example (12.1.7)]{Fulton}).  Therefore,
since all the monomial terms in the above formulae (\ref{OddchernClasses}) and
 (\ref{EvenchernClasses}) have positive
coefficients, we have that
$c_{d_1+1}(\sF^{(d_1)})\wedge\cdots\wedge c_{d_r+1}(\sF^{(d_r)})$
is not zero if
$$\left(c_2(\sF)\wedge c_1^{d_1-1}(\sF)\right)\wedge\cdots
\wedge\left(c_2(\sF)\wedge c_1^{d_r-1}(\sF)\right)$$
 is not zero.
Since the zero set of $r$ general   sections of $\sF$ vanishes
 on a subgrassmannian
$G':=\grass 2 {N+1-r}\subset G$ corresponding to an inclusion
 $\comp^{N+1-r}\subset
\comp^{N+1}$ of vector spaces, this is the same as showing that
$$c_1(\sF_{G'})^{\sum_{i=1}^r(d_i-1)}$$
is a nonzero cohomology class.
Since $\det\sF$ is the very ample line bundle on $G$ which gives
the Pl\"ucker embedding, we see that this is nonzero if
$-r+\sum_{i=1}^rd_i\le \dim
G'=2(N-1-r)$, i.e., if $\sum_{i=1}^rd_i\le 2N-2-r$.
If this is false, then since  $\sum_{i=1}^rd_i\le N$,
we conclude that
$N>2 N-2-r$ which implies that $r+1\ge N$. Since $r=N-\dim X$
this implies that $X$ is a curve.
\qed

See \cite{KS} for some related results about hypersurfaces in
$\pn 1$-bundles over
 Grassmannians.

\begin{rem*}We follow the notation used in Theorem \ref{CompIntFano} and
its proof.
If $N> \sum_{i=1}^rd_i$ then the argument of Barth and Van de Ven \cite{BVdV}
 applies to  show that there is a family of lines covering $X$ with an
$(N-\sum_{i=1}^rd_i-1)$-dimensional  space of lines
through a general point.
\end{rem*}

\section{Adjunction structure in case of $k$-very ampleness}\label{kAd}
\setcounter{theorem}{0}

Let $\widehat{X}$ be a smooth connected $n$-fold, $n\geq
3$, and let $\widehat{L}$ be a $k$-very ample line bundle on
$\widehat{X}$, $k\geq 2$. In this section we describe the first and the
second reduction of $(\widehat{X},\widehat{L})$.

We have the following general fact.

\begin{lemma}\label{nefv}Let $\widehat{X}$ be a smooth connected $n$-fold,
$n\geq
3$, and let $\widehat{L}$ be a $k$-very ample line bundle on
$\widehat{X}$, $k\geq 2$. Let $\tau$ be the nefvalue of
$(\widehat{X},\widehat{L})$. Then $\tau\leq \frac{n+1}{k}$.
\end{lemma}
\proof Let $\Phi:\widehat{X}\to W$ be the nefvalue morphism of
$(\widehat{X},\widehat{L})$. Let $C$ be an extremal curve contracted
by $\Phi$. Then $(K_\hatX+\tau \hatL)\cdot C=0$ yields
$k\tau \leq \tau \hatL\cdot C=-K_\hatX\cdot C\leq n+1$, or
$\tau\leq\frac{n+1}{k}$.\qed

We can now prove the following structure result.

\begin{theorem}\label{General}Let $\widehat{X}$ be a smooth connected
$n$-fold, $n\geq
3$, and let $\widehat{L}$ be a $k$-very ample line bundle on
$\widehat{X}$, $k\geq 2$. Then either $(\widehat{X},\widehat{L})\cong
(\pn 3,\pnsheaf 3 2)$ with $k=2$, or  the first reduction
$(X,L)$ of $(\widehat{X},\widehat{L})$ exists    and $\widehat{X}\cong X$,
$L\cong \hatL$.
Furthermore either:
\begin{enumerate}
\item[{\rm i)}] $(\widehat{X},\widehat{L})\cong (\pn 3,\sO_{\pn
3}(2))$, $k=2$;
\item[{\rm ii)}] $(\widehat{X},\widehat{L})\cong (\pn 3,\sO_{\pn
3}(3))$, $k=3$;
\item[{\rm iii)}] $(\widehat{X},\widehat{L})\cong (\pn 4,\sO_{\pn
4}(2))$, $k=2$;
\item[{\rm iv)}] $(\widehat{X},\widehat{L})\cong (Q,\sO_Q(2))$, $Q$
hyperquadric in $\pn 4$, $k=2$;
\item[{\rm v)}] there exists a morphism $\psi:\hatX\to C$ onto a smooth
curve $C$ such that $2K_\hatX+3\hatL\approx \psi^*H$ for some ample line bundle
$H$ on $C$ and $(F,{\hatL}_F)\cong(\pn 2,\sO_{\pn 2}(2))$ for any fiber
$F$ of $\psi$, $k=2$;
 \item[{\rm vi)}] $(\widehat{X},\widehat{L})$ is
a Mukai variety, i.e., $K_\hatX\approx -(n-2)\hatL$ and either $n=4,5$ and $k=2$
or $n=3$ and $k\leq 4$;
 \item[{\rm vii)}] $(\widehat{X},\widehat{L})$
is a Del Pezzo fibration over a curve such that $(F,{\widehat{L}}_F)\cong
(\pn 3,\sO_{\pn 3}(2))$ for each general fiber $F$, $n=4$, $k=2$,
\end{enumerate}
or there exists the second reduction $(Z,\sD)$, $\vphi:X\to Z$ of
 $(\widehat{X},\widehat{L})$. In this case the following hold.
\begin{enumerate}
\item[{\rm 1)}] If $n\geq 4$, then $X\cong Z$;
\item[{\rm 2)}] If $n=3$, then either $X\cong Z$ or $k=2$ and $\vphi$
only contracts divisors $D\cong \pn 2$ such that $L_D\cong \sO_{\pn
2}(2)$; furthermore $\sO_D(D)\cong \sO_D(-1)$ and $Z$ is smooth.
\end{enumerate}\end{theorem}
\proof We use general results from adjunction theory for which we refer
to \cite{Book}. From \cite[(9.2.2)]{Book} we know that
$K_{\widehat{X}}+(n-1)\widehat{L}$ is spanned unless either
$(\widehat{X},\widehat{L})\cong (\pn n,\sO_{\pn n}(1))$, or
$\widehat{X}\subset \pn {n+1}$ is a quadric hypersurface and
$L\approx\sO_{\pn {n+1}}(1)_{|\widehat{X}}$, or
$(\widehat{X},\widehat{L})$ is  a scroll over a curve. Since $L\cdot C\geq 2$
for each curve $C$ on $\widehat{X}$, all the above cases are excluded.

Thus we can conclude that $K_{\widehat{X}}+(n-1)\widehat{L}$ is spanned.
Then from \cite[(7.3.2)]{Book} we know that $K_{\widehat{X}}+(n-1)
\widehat{L}$ is big unless either $(\widehat{X},\widehat{L})$ is a Del
Pezzo variety, i.e., $K_{\widehat{X}}\approx -(n-1)\widehat{L}$, or
$(\widehat{X},\widehat{L})$ is a quadric fibration over a smooth curve, or
$(\widehat{X},\widehat{L})$ is a scroll over a normal surface. Then, as
above, the quadric fibration and the scroll cases are excluded, so that
$(\widehat{X},\widehat{L})$ is a Del Pezzo variety.
In this case $\tau=n-1$, so that Lemma (\ref{nefv}) gives $2\leq k\leq
\frac{n+1}{n-1}$. Hence $n=3$. By looking over the list of Del Pezzo
$3$-folds (see \cite[(8.11)]{Fujita}) we conclude that
$(\widehat{X},\widehat{L})\cong (\pn 3,\sO_{\pn 3}(2))$ in this case.

Thus we can assume that the first reduction, $(X,L)$, of
$(\widehat{X},\widehat{L})$ exists and in fact $\widehat{X}\cong X$, since
otherwise we can find a line $\ell$ on $\widehat{X}$ such that
$\widehat{L}\cdot \ell=1$. From \cite[(7.3.4), (7.3.5), (7.5.3)]{Book}) we know
that on $\widehat{X}\cong X$ the line bundle $K_X+(n-2)L$ is nef and big unless
either \begin{enumerate}
\item[{\rm a)}] $(X,L)\cong (\pn 3,\sO_{\pn
3}(3))$;
\item[{\rm b)}] $(X,L)\cong (\pn 4,\sO_{\pn
4}(2))$;
\item[{\rm c)}] $(X,L)\cong (Q,\sO_Q(2))$, $Q$
hyperquadric in $\pn 4$;
\item[{\rm d)}] there exists a morphism $\psi:X\to C$ onto a smooth
curve $C$ such that $2K_X+3L\approx \psi^*H$ for some ample line bundle
$H$ on $C$ and $(F,L_F)\cong(\pn 2,\sO_{\pn 2}(2))$ for any fiber
$F$ of $\psi$;
\item[{\rm e)}] $K_X\approx -(n-2)L$, i.e., $(X,L)$ is a Mukai variety;
\item[{\rm f)}] $(X,L)$ is a Del Pezzo fibration over a smooth curve under
the morphism, $\Phi_L$, given by $|m(K_X+(n-2)L)|$ for $m\gg 0$;
\item[{\rm g)}] $(X,L)$ is a quadric fibration over a normal surface under
$\Phi_L$; or
\item[{\rm h)}] $(X,L)$ is a scroll over a normal threefold under
$\Phi_L$.\end{enumerate}
Cases a), b), c), d), e) lead to cases ii), iii), iv), v), vi)
respectively.

Case f) gives case vii). To see this, let $F$ be a general fiber of $\Phi_L$
and let $L_F$ be the restriction of $L$ to $F$. Let $\tau_F$ be the nefvalue
of $(F,L_F)$. Then $K_F+(n-2)L_F$ is trivial, and hence $\tau_F=n-2=\dim
F-1$. Therefore the same argument as above, by using again
\cite[(8.11)]{Fujita}, gives $\dim F=3$ and $(F,L_F)\cong (\pn 3,\sO_{\pn
3}(2))$.

Note that in case e), one has
$\tau=n-2$, so that Lemma (\ref{nefv}) yields $2\leq k\leq \frac{n+1}{n-2}$.
Thus either $n=4,5$ and $k=2$, or $n=3$ and $k\leq 4$.

In cases g), h) we can find a line $\ell$ on $X$ such that $L\cdot \ell=1$, so
that they are excluded.

Thus we can assume that the second reduction, $(Z,\sD)$, $\vphi:X\to Z$, of
$(\widehat{X}, \widehat{L})$ exists. Use the structure results of the
second reduction (see \cite[(7.5.3), (12.2.1)]{Book}). If $n\geq 4$ we see
that we can always find a line  $\ell$ on $X$ such that $L\cdot\ell=1$.
Then $X\cong Z$. If $n=3$, either $X\cong Z$ or $\vphi$ contracts divisors
$D\cong \pn 2$ such that $L_D\cong \sO_{\pn 2}(2)$. Then $\sO_D(D)\cong
\sO_D(-1)$ and $Z$ is smooth in this case.\qed

The following is an immediate consequence of Theorem (\ref{General}).

\begin{corollary} \label{ample}Let $\widehat{X}$ be a smooth connected
$n$-fold, $n\geq
3$, and let $\widehat{L}$ be a $k$-very ample line bundle on
$\widehat{X}$, $k\geq 2$. Then $K_{\widehat{X}}+(n-2)\widehat{L}$ is ample
if $n\geq 4$ unless $k=2$ and either $(\widehat{X},\widehat{L})\cong (\pn
4,\sO_{\pn 4}(2))$, or $n=4,5$ and $(\widehat{X},\widehat{L})$ is a Mukai
variety, or $n=4$ and $(\widehat{X},\widehat{L})$ is a Del Pezzo fibration
over a curve such that $(F, {\widehat{L}}_F)\cong(\pn 3,\sO_{\pn 3}(2))$ for
each fiber $F$.\end{corollary}

The results of this section justify the study of Mukai varieties of
dimension $n=3,4,5$, polarized by a $k$-very ample line bundle, $k\geq 2$.

\section{Mukai varieties of dimension $n\geq 4$}\label{Mukai}
\setcounter{theorem}{0}

Let $(X,L)$ be a Mukai variety of dimension $n\geq 3$, i.e.,
$K_X\approx -(n-2)L$, polarized by a $k$-very ample line bundle $L$, $k\geq
2$ (see  \cite{Mu1}, \cite{Mu2} for classification results of Mukai
varieties).
Since the nefvalue, $\tau$, of such pairs $(X,L)$ is $\tau=n-2$, we
immediately see from Lemma (\ref{nefv}) that either $n=4,5$, $k=2$, or
$n=3$, $k\leq 4$ (compare with the proof of (\ref{General})). We have the
following
result.

\begin{theorem}\label{MuThm} Let $(X,L)$ be a Mukai variety of dimension
$n\geq 4$ polarized by a $k$-very ample line bundle $L$, $k\geq 2$. Then
either
\begin{enumerate}
\em\item\em $n=4$, $k=2$,  $(X,L)\cong (Q,\sO_Q(2))$, $Q$ hyperquadric in
$\pn 5$, or
\em\item\em $n=5$, $k=2$,  $(X,L)\cong (\pn 5, \sO_{\pn
5}(2))$.\end{enumerate}\end{theorem}
\proof By the above we know that $k=2$ and $n=4,5$. Let $V$ be the $3$-fold
section obtained as transversal intersection of $n-3$ general members of
$|L|$. Let $L_V$ be the restriction of $L$ to $V$. Note that $K_V\approx
-L_V$, so that $V$ is a Fano $3$-fold, and $L_V$ is $k$-very ample.

We denote by $r$ the index of $V$. Then $L_V\approx -K_V=rH$ for some ample
line bundle $H$ on $V$. Note that we have $r\geq 2$ since otherwise Shokurov's
theorem \cite{Sho} (see also \cite{R}) applies to say that either $V=\pn 1\times
\pn 2$ or there exists a line $\ell$ on $V$ with respect to $H$. In the latter
case $(H\cdot\ell)_V=(L\cdot\ell)_X=1$, which contradicts the assumption $k\geq
2$. The following argument rules out the former case $V=\pn 1\times \pn 2$.

If $V=\pn 1\times \pn 2$, as corollary of the extension theorem \cite[Prop.
III]{SoAmple} we conclude that $X$ is a linear $\pn 3$-bundle over $\pn 1$,
$p:X\to \pn 1$, with the restriction $p_V$ giving the map  $V=\pn
1\times\pn 2\to \pn 1$. By taking the direct image of
$$0\to \sO_X\to L\to L_V\to 0$$
we get the exact sequence
$$0\to \sO_{\pn 1}\to \sE\to \sE/\sO_{\pn 1}\to 0,$$
where $\sE:=p_*L$. Since ${\Bbb P}(\sE/\sO_{\pn 1})\cong \pn 1\times \pn 2$
we conclude that $\sE/\sO_{\pn 1}=\sO_{\pn 1}(1)\oplus\sO_{\pn 1}(1)\oplus
\sO_{\pn 1}(1)$. Thus $\deg(\sE/\sO_{\pn 1})=3=\deg(\sE)< {\rm rank}(\sE)=4$.
Therefore $L$ cannot be ample.

By Lefschetz theorem we have $H\approx H'_V$ for some line bundle $H'$ on $X$,
as well as $L\approx rH'$. Hence in particular $H'$ is ample. We have
$K_X+r(n-2)H'\approx \sO_X$, so that $r(n-2)\leq n+1$ by a well known result
due to Maeda  (see e.g., \cite[(7.2.1)]{Book}). If $r=4,3$ we find numerical
contradictions since we are assuming $n\geq 4$. Thus $r=2$ and $n\leq 5$.

If $n=5$ we have $K_X\approx -6H'$ and if $n=4$ we have $K_X\approx -4H'$.
By the Kobayashi-Ochiai theorem (see e.g., \cite[(3.6.10)]{Book}) we get in
the former case $(X,H')\cong (\pn 5,\sO_{\pn 5}(1))$, or $(X,L)\cong
(\pn 5,\sO_{\pn 5}(2))$ as in case $1$) of the statement, and in the latter
case $(X,H')\cong (Q,\sO_Q(1))$, $Q$ hyperquadric in $\pn 5$, or
$(X,L)\cong (Q,\sO_Q(2))$ as in case $2$) of the statement.
\qed

\begin{rem*}\label{j1} Note that in both cases 1), 2) of Theorem
(\ref{MuThm}) the line bundle $L$ is in fact $2$-jet ample (see
\cite[Corollary (2.1)]{Plenum}).
\end{rem*}

\section{The Fano $3$-fold case}\label{FanoS}
\setcounter{theorem}{0}

In this section we classify the $3$-dimensional  Mukai varieties $(X,L)$
polarized by a $k$-very ample line bundle $L$, $k\geq 2$, i.e., we classify
all Fano $3$-folds $X$ such that the anticanonical divisor $-K_X$ is
$k$-very ample, $k\geq 2$.
\begin{prgrph*}{A special case}\label{SpecialEx} Let us start by studying a
particular case. This case has a special interest also because it gives a
simple explicit example of a line bundle which is $2$-very ample but
not $2$-jet ample. This example is case 4) in the Iskovskih-Shokurov's list
\cite[Table 21]{IS} of Fano $3$-folds of first species, i.e., $b_2(X)=1$.

\begin{proposition}\label{SpecialProp} Let $X$ be a smooth double cover of $\pn
3$,
$p:X\to \pn 3$, branched along a quartic. Then $L:=-K_X$ is $2$-very ample
but not $2$-jet ample.
\end{proposition}
\proof We have $L:=-K_X=p^*\sO_{\pn 3}(2)$. First we show that $L$ is not
$2$-jet ample. Let $R$ be the ramification divisor of $p$. Let $\sZ$ be a
length $3$ zero dimensional subscheme of $X$ such that ${\rm
Supp}(\sZ)=\{x\}$ with $x\in R$. We can assume that $R$ is defined at $x$
by a local coordinate, $s$, i.e., $R=\{s=0\}$ at $x$. Consider local
coordinates $(s,v,w)$ on $X$ at $x$. Let $y\in \pn 3$ be a point, belonging
to the branch locus of $p$, such that $y=p(x)$. We can consider local
coordinates $(t,v,w)$ on $\pn 3$ at $y$, where $p^*t=s^2$. We have
\begin{equation}\label{Decomp}
H^0(L)\cong H^0(p_*p^*\sO_{\pn 3}(2))\cong H^0(\sO_{\pn 3}(2))\oplus
H^0(\sO_{\pn 3}).
\end{equation}
Therefore we can find a base, $\sB$, of $H^0(L)$ given by the pullback of
sections of $\sO_{\pn 3}(2)$ and one more section, $\grs\in
H^0(\sO_{\pn 3})$, which in local coordinates around $x$ is of the form
 $\lambda s$ with $\lambda$ a holomorphic function that doesn't vanish at $x$.
That is,
recalling that $s^2=p^*t$,
$v=p^*v$,
$w=p^*w$,
$$\sB=<1,s^2,v,w,s^4,v^2,w^2,s^2v,s^2w,vw,\grs>.$$
On the other hand, $H^0(L/{\frak m}_x^3)(=H^0(\sO_\sZ(L)))$ contains the
elements $sv$, $sw$ which are not images of elements of the base $\sB$.
This shows that the restriction map $H^0(L)\to H^0(\sO_\sZ(L))$ is not
surjective. Thus $L$ is not $2$-jet ample.

We prove now that $L$ is $2$-very ample. Consider a $0$-dimensional
subscheme $\sZ$ of $X$ of length $3$. Recalling (\ref{Decomp}), the fact
that $\sO_{\pn 3}(2)$ is $2$-very ample and Lemma (\ref{GenFa}), we see
that the restriction map $H^0(L)\to H^0(\sO_\sZ(L))$ is always surjective
except possibly in the case when ${\rm Supp}(\sZ)$ is a single point, $x$,
belonging to the ramification divisor of the cyclic covering $p$.

Thus, let us assume ${\rm Supp}(\sZ)=\{x\}$ and consider the ideals
$\sJ_i:=(\sJ_\sZ,{\frak m}_x^i)$, where $\sJ_\sZ$, ${\frak m}_x$ denote
the ideal sheaves of $\sZ$, $x$ respectively, the sheaves
$\sO_i:=\sO_X/\sJ_i$, the maps
$p_i:\sO_\sZ\to \sO_i$, $i=1,2,3$, and the cofiltration
$\sO_3\to\sO_2\to\sO_1\to 0.$
Note that the following hold true.
\begin{itemize}
\item If $H^0(\sO_\sZ)$ is generated by only terms of degree $\leq 1$ in
the local
coordinates $s$, $v$, $w$ on $X$ at the point $x$ then (as observed before
in the proof
that $L$ is not $2$-jet ample) the image of $H^0(L)$ can generate
$H^0(\sO_\sZ(L))$, so
the restriction map $H^0(L)\to H^0(\sO_\sZ(L))$ is surjective in this case
and we are
done;
\item ${\rm length}(\sO_1)=1$;
\item ${\rm length}(\sO_i)\neq{\rm length}(\sO_{i+1})$, $i=1,2$. Indeed,
otherwise, $(\sJ_\sZ,{\frak m}_x^i)=(\sJ_\sZ,{\frak m}_x^{i+1})$ so
that ${\frak m}_x^i\subset \sJ_{\sZ}$ and therefore $H^0(\sO_\sZ)$ is
generated by only constant terms or linear terms. By the above we are done in
this case.\end{itemize} Thus we are reduced to consider the case when ${\rm
length}(\sO_2)=2$, ${\rm length}(\sO_3)=3$.

We claim that $H^0(\sO_\sZ)$ contains at least one quadratic term in $s$,
$v$, $w$.
Indeed, if not, ${\frak m}_x^2\subset \sJ_\sZ$ and hence we would have
$(\sJ_\sZ,
{\frak m}_x^2)=(\sJ_\sZ,{\frak m}_x^3)$, which gives the contradiction
${\rm length}(\sO_2)={\rm length}(\sO_3)$.

Furthermore, since ${\rm length}(\sO_X/(\sJ_\sZ,{\frak m}_x^2))=2$ and
${\rm length}(\sO_X/{\frak m}_x^2)=4$, we conclude that $\sJ_\sZ$ contains two
independent linear terms, say $f$, $g$, not belonging to ${\frak m}_x^2$. Write
$$f=as+bv+cw,\;\;g=ds+ev+hw,$$
where the coefficients $a$, $b$, $c$, $d$, $e$, $h$ belong to $\sO_{X,x}$.
Let $\sB$ be
the base of $H^0(L)$ constructed in the first part of the proof, where we
showed that $L$
is not $2$-jet ample. Following that argument we see that $L$ is $2$-very
ample as soon as we show that the elements $sv$, $sw$ can be written in
$\sO_\sZ$ in terms of elements of $\sB$.

We go on by a case by case analysis. First, assume $a\neq 0$, i.e., $a$
invertible in
$\sO_{X,x}$ and write
$$asw=w(as+bv+cw)-bvw-cw^2.$$
Since $as+bv+cw=f=0$ in $\sO_\sZ$, up to dividing by $a$, we can express
$sw$ in terms of
$vw, w^2\in \sB$ in $\sO_\sZ$. Similarly, writing
$$asv=v(as+bv+cw)-bv^2-cvw,$$
we conclude that $sv$ can be expressed in terms of $v^2, vw\in\sB$ in
$\sO_\sZ$.
If $d\neq 0$ we get the same conclusion.

Thus it remains to consider the case when $a=d=0$. In this case $f=bv+cw$,
$g=ev+hw$ in
$\sO_\sZ$. Then, solving with respect to $v$, $w$, and noting that
$\left(\begin{array}{ll} b &c\\ e &h\end{array}\right)$ is a rank two
matrix since $f$,
$g$ are independent, we can express $v$, $w$ as linear functions of $f$,
$g$ in $\sO_\sZ$.
Since $f, g\in \sJ_\sZ$, we conclude that $v, w\in \sJ_\sZ$ and hence $sv$,
$sw$ belong to
$\sJ_\sZ$. Therefore $sv=sw=0$ in $\sO_\sZ$. \qed
\end{prgrph*}

The following general result is a  consequence of Fujita's classification
\cite{Fupapers}, \cite[(8.11)]{Fujita}  of Del Pezzo $3$-folds (see also
\cite{MoMu}, \cite{IS} and \cite{Murre} for a complete classification of Fano
$3$-folds).

Note that in each case of the theorem below the line bundle $L$ is in fact
$k$-very ample (see \cite[Corollary (2.1)]{Plenum}, Lemma (\ref{Tensor}) and
Proposition (\ref{SpecialProp}).

\begin{theorem}\label{Fano} Let $X$ be a Fano threefold. Assume that $L:=-K_X$
is $k$-very ample, $k\ge 2$. Then either:
\begin{enumerate}
\em\item\em $X$ is a divisor on $\pn 2\times\pn 2$ of bidegree $(1,1)$,
$L=\sO_X(2,2)$,
$k=2$;
\em\item\em $X=\pn 1\times\pn 2$, $L=\sO_{\pn 1}(2)\bx\sO_{\pn 2}(3)$, $k=2$;
\em\item\em $X=V_7$, the blowing up of $\pn 3$ at a point, $L=2(q^*\sO_{\pn
3}(2)-E)$,
$q:V_7\to \pn 3$, $E$ the exceptional divisor, $k=2$;
\em\item\em $X=\pn 1\times\pn 1\times\pn 1$, $L=\sO_{\pn 1}(2)\bx\sO_{\pn
1}(2)\bx\sO_{\pn
1}(2)$, $k=2$;
\em\item\em $X=\pn 3$, $L=\sO_{\pn 3}(4)$, $k=4$;
\em\item\em $X$ is a hyperquadric in $\pn 4$, $L=\sO_X(3)$, $k=3$;
\em\item\em $X$ is a cubic hypersurface in $\pn 4$, $L=\sO_X(2)$, $k=2$;
\em\item\em $X$ is the complete intersection of two quadrics in $\pn 5$,
$L=\sO_X(2)$,
$k=2$;
\em\item\em $X$ is a double cover of $\pn 3$, $p:X\to \pn 3$, branched
along a quartic;
$L=p^*\sO_{\pn 3}(2)$ is $2$-very ample but not $2$-jet ample;
\em\item\em $X$ is the section of the Grassmannian $\grass 2 5$ {\rm (}of
lines in $\pn 4${\rm )} by a linear subspace of codimension $3$,
$L=\sO_X(2)$, $k=2$.
\end{enumerate}
\end{theorem}
 \proof Let $r$ be the
index of $X$.  Then $L:=-K_X=rH$ for some ample line bundle $H$ on $X$. Note
that we have $r\geq 2$ since otherwise Shokurov's theorem \cite{Sho}
applies to say that either $X=\pn 1\times\pn 2$
 or there exists a line $\ell$ with respect to $H$. In the former case we
are in case
$2$) of the statement. In the latter case
$H\cdot\ell=L\cdot\ell=1$, which contradicts the assumption $k\geq 2$.

If $r=4,3$, by using the Kobayashi-Ochiai theorem (see e.g.,
\cite[(3.6.1)]{Book}) we find  cases $5$), $6$) of the statement
respectively.

Thus we can assume $r=2$. In this case $(X,H)$ is a Del Pezzo $3$-fold
described as in  \cite[(8.11)]{Fujita}. A direct check shows that the cases
listed in \cite[(8.11)]{Fujita}
 lead to cases $1$), $3$), $4$), $7$), $8$), $9$), $10$) of the statement.
Recall that case $9$) is discussed in Proposition (\ref{SpecialProp}).

Notice that the case of $X={\Bbb P}(\sT)$, for the tangent bundle $\sT$ of $\pn
2$, as in \cite[(8.11), $6$)]{Fujita} gives our case $1$) (see Remark (\ref{PT})
below). Note also that case $1$) of \cite[(8.11)]{Fujita}, when $(X,H)$ is a
weighted hypersurface of degree $6$ in the weighted projective space ${\Bbb
P}(3,2,1,\ldots,1)$ with $H^3=1$, is ruled out since $L=2H$ is not even very
ample. To see this notice that there exist a smooth surface $S$ in $|H|$ and a
smooth curve $C$ in $|H_S|$ (see \cite[(6.1.3), (6.14)]{Fujita}). On $S$ we
have $K_S\approx -H_S$, so that $K_S^2=1$. Therefore $C$ is a smooth elliptic
curve with $H\cdot C=H^3=1$, i.e., $L\cdot C=2$.\qed

\begin{rem*}\label{PT} It is a standard fact that ${\Bbb P}(\sT)$, for the
tangent bundle $\sT:=\sT_{\pn n}$ of $\pn n$, is embedded in $\pn n\times \pn n$ as a 
divisor of bidegree $(1,1)$ (see also \cite{Sato}). To see this, note that $\sT(-1)$ 
is spanned with $n+1$ sections. Thus letting $\xi$ denote the tautological line bundle 
of $\proj{\sT(-1)}$, the map  $f:{\Bbb  P}(\sT)\to \pn n$ associated to $|\xi|$ is an 
embedding on fibers of the bundle projection $p:{\Bbb P}(\sT)\to \pn n$. The product map 
$(f,p):{\Bbb P}(\sT)\times {\Bbb P}(\sT)\to \pn n\times \pn n$ is thus an embedding with 
image a divisor $D$ such that $\sO_{\pn n\times \pn n}(D)_{|F}\cong\pnsheaf {n-1}1$ on the 
fibers $F$ of $p$. The fibers $(F',\xi_{F'})$ of $f$ are $\cong \pnpair {n-1}1$. To see this 
for $F'$, a generic fiber of $f$, note that 
$c_1(p^*\pnsheaf n 1)^{n-1}\cdot F'=c_1(p^*\pnsheaf n 1)^{n-1}\cdot c_1(\xi)^n$.  Using the
defining equation for the Chern classes of $\sT(-1)$, we see that this equals
$p^*\left(c_1(\pnsheaf n 1)^{n-1}\cdot c_1(\sT(-1))\right)\cdot c_1(\xi)^{n-1}=
c_1(\pnsheaf n 1)^{n-1}\cdot c_1(\sT(-1))=1.$ Since $D$, $f(D)$, and the generic fiber of 
$f$  are all connected, it follows that all fibers of $f$ are connected. From this it follows 
that all fibers of $f$ are isomorphic if the automorphism group of $D$ acts transitively on 
$D$. This can be seen by observing that given two nonzero tangent vectors of $\pn n$ there is 
an automorphism of $\pn n$ which takes one tangent vector to the other.  \end{rem*}

\begin{rem*}\label{kspanj} Note that the line bundles $L$ in (\ref{MuThm})
and (\ref{Fano}) are of course $k$-spanned and also $k$-jet ample, with the
only exception of Case $9$) in (\ref{Fano}) (see (\ref{kthemb}),
(\ref{SpecialProp}) and \cite[Corollary (2.1)]{Plenum}). Thus
 (\ref{MuThm}) and (\ref{Fano}) also give the
classification of Mukai varieties $(X,L)$ of dimension $n\geq 3$ polarized
by either a $k$-spanned or a $k$-jet ample line bundle $L$, with the only
exception, for $k$-jet ampleness, of Case $9$) in (\ref{Fano}).
\end{rem*}

\bigskip

{
\begin{tabular}{l} Mauro C. Beltrametti\\
 Dipartimento di Matematica\\
Via Dodecaneso 35\\
 I-16146 Genova, Italy\\
beltrame@dima.unige.it\\
\ \ \ \\
Sandra Di Rocco\\
Department of Mathematics\\
KTH, Royal Institute of Technology\\
S-100 44 Stockolm, Sweden\\
sandra@math.kth.se\\
\ \ \ \\
Andrew J. Sommese\\
Department of Mathematics\\
Notre Dame, Indiana, 46556, U.S.A\\
sommese@nd.edu\\
http://www.nd.edu/$\sim$sommese/index.html
 \end{tabular}
}


\begin{thebibliography}{999}
\bibitem{BVdV} W. Barth  and A. Van de Ven, ``Fano-varieties of lines on
hypersurfaces,''
Arch.\ Math.\ (Basel)  31
(1978), 96--104.

\bibitem{Duke}  M.C. Beltrametti, P. Francia, and A.J. Sommese, ``On
Reider's method and higher order embeddings,'' Duke Math.\ J. 58 (1989),
425--439.

\bibitem{BSS2} M.C. Beltrametti, M. Schneider, and A.J. Sommese,
``Threefolds of degree $11$ in $\pn 5$,''
in {\em Complex Projective Geometry}, ed.\ by G. Ellingsrud, C. Peskine, G.
Sacchiero, and S.A. Stromme,
London Math.\ Soc. Lecture Note Ser.
179 (1992), 59--80.

\bibitem{BSAq}  M.C. Beltrametti and A.J. Sommese, ``On $k$-spannedness
for projective
surfaces,'' in {\em Algebraic Geometry, Proceedings of Conference on
Hyperplane Sections,
L'Aquila, Italy, 1988}, ed.\ by A.J. Sommese, A. Biancofiore, and E.L.
Livorni,  Lecture
Notes in Math. 1417 (1990),  24--51, Springer-Verlag, New York.

\bibitem{BSCo}  M.C. Beltrametti and A.J. Sommese, ``Zero cycles and $k$-th
order embeddings of smooth projective surfaces (with an appendix by L.
G\"ottsche),''
in {\em 1988 Cortona Proceedings on Projective Surfaces and their
Classification}, ed.\ by F. Catanese and C. Ciliberto, Sympos.\ Math. 32
(1992), 33--48,
INDAM, Academic Press, London.

\bibitem{BSMZ} M.C. Beltrametti and A.J. Sommese, ``On the preservation of
$k$-very ampleness under adjunction,'' Math.\ Z. 212 (1993), 257--283.

\bibitem{Plenum} M.C. Beltrametti and A.J. Sommese, ``On $k$-jet
ampleness,'' {\em Complex Analysis and Geometry}, ed. by V. Ancona and A.
Silva, The University Series in Mathematics, Plenum Press, (1993), 355--376.

\bibitem{Book} M.C. Beltrametti and A.J. Sommese, ``{\em The Adjunction
Theory of Complex Projective Varieties}, Expositions in Mathematics, 16, Walter
de Gruyter, Berlin (1995).

\bibitem{Embedding} M.C. Beltrametti and A.J. Sommese, ``Notes on embeddings
of blowups,''   Journal of Algebra 186 (1996), 861--871.

\bibitem{CG} F. Catanese and L. G\"ottsche, ``$d$-very-ample line bundles and
embeddings of Hilbert schemes of $0$-cycles,'' Manuscripta   Math. 68
(1990),  337--341.

\bibitem{Sandra1} S. Di Rocco, ``$k$-very ample line bundles on Del Pezzo
surfaces,''  Math.\ Nachr. 179 (1996), 47--56.

\bibitem{Fulton} W. Fulton, {\em Intersection Theory}, Ergeb.\ Math.\
Grenzgeb.\ (3)
 2, Springer-Verlag, Berlin, (1984).

\bibitem{Fupapers} T. Fujita, ``On the structure of polarized manifolds with
total deficiency one, I; II; III,'' J. Math. Soc. Japan 32 (1980), 709--725;
33 (1981), 415--434; 36 (1984), 75--89.

\bibitem{Fujita} T. Fujita, {\em
Classification Theories of Polarized Varieties}, London Math.\ Soc.\ Lecture
Note Ser. 155, Cambridge University Press, (1990).

\bibitem{IS} V.A. Iskovskih and V.V. Shokurov,  ``Biregular theory of
Fano $3$-folds,'' in {\em  Proceedings of the Algebraic Geometry Conference,
Copenhagen 1978},   Lecture Notes in Math. 732 (1979),  171--182,
Springer-Verlag, New York.


\bibitem{KS}M. Kim and A.J. Sommese, ``Two results on
branched coverings of Grassmannians,''  preprint, 1996.

\bibitem{MoMu} S. Mori and S. Mukai, ``Classification of Fano $3$-folds with
$B_2 \ge 2$,''  Manuscripta   Math. 36 (1981),  147--162.

\bibitem{Mu1} S. Mukai, ``New classification of Fano threefolds and
manifolds of coindex $3$,''
preprint, 1988.

\bibitem{Mu2} S. Mukai,``Biregular classification of Fano threefolds and
Fano manifolds of coindex $3$,'' Proc.\ Nat.\ Acad.\ Sci.\ U.S.A. 86
(1989),  3000--3002.

\bibitem{Murre} J.P. Murre, ``Classification of Fano threefolds according to
Fano and Iskovskih,'' in {\em Algebraic $3$-folds, Proceedings Varenna, 1981},
ed.\ by A. Conte,  Lecture Notes in Math. 947 (1982),  35--92,
Springer-Verlag, New York.

\bibitem{R} M. Reid, ``Lines on Fano $3$-folds according to Shokurov,''
Mittag-Leffler Report n.\ 11 (1980).

\bibitem{Sato} E. Sato, ``Varieties which have two projective bundle
structures,'' J. Math. Kyoto 25 (1985), 445--457.

 \bibitem{Sho} V. Shokurov, ``The existence of a
straight line on Fano $3$-folds,'' Izv. Akad. Nauk. 43 (1979), engl. trans.
Math. U.S.S.R. Izv. 15 (1980).

\bibitem{SoAmple} A.J. Sommese, ``On manifolds that cannot be ample divisors,''
Math. Ann. 221 (1976), 55--72.

\end{thebibliography}
\end{document}